\documentclass[a4paper,12pt]{amsart}

\usepackage{amsmath}
\usepackage{amssymb}
\usepackage{mathrsfs}
\usepackage{ifthen}
\usepackage{graphicx}
\usepackage[T1]{fontenc} 

\def\switchlinenumbers{\@ifstar
    {\let\makeLineNumberOdd\makeLineNumberRight
     \let\makeLineNumberEven\makeLineNumberLeft}%
    {\let\makeLineNumberOdd\makeLineNumberLeft
     \let\makeLineNumberEven\makeLineNumberRight}%
    }

\def\setmakelinenumbers#1{\@ifstar
  {\let\makeLineNumberRunning#1%
   \let\makeLineNumberOdd#1%
   \let\makeLineNumberEven#1}%
  {\ifx\c@linenumber\c@runninglinenumber
      \let\makeLineNumberRunning#1%
   \else
      \let\makeLineNumberOdd#1%
      \let\makeLineNumberEven#1%
   \fi}%
  }

\nonstopmode \numberwithin{equation}{section}
\newtheorem*{theorem*}{Theorem}

\newtheorem{thm}{Theorem}[section]
\newtheorem{cor}[equation]{Corollary}
\newtheorem{lem}[equation]{Lemma}
\newtheorem{prop}[equation]{Proposition}

\theoremstyle{definition}
\newtheorem{defn}{Definition}[section]

\newtheorem{prob}[equation]{Problem}
\newtheorem{rem}{Remark}[section]

\newcounter{minutes}
\divide\time by 60
\newcounter{hours}
\multiply\time by 60
\addtocounter{minutes}{-\time}

\newcounter {own}
\def\theown {\thesection       .\arabic{own}}

\newenvironment{pf}[1][]{%
 \vskip 3mm
 \noindent
 \ifthenelse{\equal{#1}{}}%
  {{\slshape Proof. }}%
  {{\slshape #1.} }%
 }%
{\qed\bigskip}

\newcounter{alphabet}

\def\be{\begin{equation}}
\def\ee{\end{equation}}

\newcommand{\bee}{\begin{enumerate}}
\newcommand{\eee}{\end{enumerate}}

\newcommand{\blem}{\begin{lem}}
\newcommand{\elem}{\end{lem}}
\newcommand{\bthm}{\begin{thm}}
\newcommand{\ethm}{\end{thm}}
\newcommand{\bcor}{\begin{cor}}
\newcommand{\ecor}{\end{cor}}
\newcommand{\beg}{\begin{examp}}
\newcommand{\eeg}{\end{examp}}
\newcommand{\begs}{\begin{examples}}
\newcommand{\eegs}{\end{examples}}
\newcommand{\bdefe}{\begin{defin}}
\newcommand{\edefe}{\end{defin}}
\newcommand{\bprob}{\begin{prob}}
\newcommand{\eprob}{\end{prob}}
\newcommand{\bei}{\begin{itemize}}
\newcommand{\eei}{\end{itemize}}

\newcommand{\norm}[1]{\left\lVert#1\right\rVert}

\begin{document}

\title{Arithmetic Bohr radius for the Minkowski space}

\author{Vasudevarao Allu}
\address{Vasudevarao Allu,
	Department of Mathematics,
School of Basic Sciences,
Indian Institute of Technology Bhubaneswar,
Bhubaneswar-752050, Odisha, India.}
\email{avrao@iitbbs.ac.in}

\author{Himadri Halder}
\address{Himadri Halder,
Department of Mathematics,
Indian Institute of Technology Bombay, Powai, Mumbai,
Maharashtra-400076, India.}
\email{himadrihalder119@gmail.com, himadri@math.iitb.ac.in}

\author{Subhadip Pal}
\address{Subhadip Pal,
	Department of Mathematics,
	School of Basic Sciences,
	Indian Institute of Technology Bhubaneswar,
	Bhubaneswar-752050, Odisha, India.}
\email{subhadippal33@gmail.com}

\subjclass[{AMS} Subject Classification:]{Primary 32A05, 32A10, 32A17; Secondary 30B10}
\keywords{Bohr radius, Arithmetic Bohr radius, Holomorphic functions, Reinhardt domain}

\def\thefootnote{}
\footnotetext{ {\tiny File:~\jobname.tex,
printed: \number\year-\number\month-\number\day,
          \thehours.\ifnum\theminutes<10{0}\fi\theminutes }
} \makeatletter\def\thefootnote{\@arabic\c@footnote}\makeatother

\begin{abstract}
	The main aim of this paper is to study the arithmetic Bohr radius for holomophic functions defined on a Reinhardt domain in $\mathbb{C}^n$ with positive real part. The present investigation is motivated by the work of Lev Aizenberg [Proc. Amer. Math. Soc. 128 (2000), 2611--2619]. A part of our study in the present paper includes a connection between the classical Bohr radius and the arithmetic Bohr radius of unit ball in the Minkowski space $\ell^n_{q}\, , 1\leq q\leq \infty$. Further, we determine the exact value of a Bohr radius in terms of arithmetric Bohr radius.

\end{abstract}

\maketitle
\pagestyle{myheadings}
\markboth{Vasudevarao Allu, Himadri Halder, and Subhadip Pal}{Arithmetic Bohr radius for the Minkowski space}

\section{Introduction}

A domain $\Omega$ centered at the origin in $\mathbb{C}^n$ is said to be a complete Reinhardt domain if for all $z=(z_{1},\ldots,z_{n}) \in \Omega$ and $\xi_{i} \in \overline{\mathbb{D}}$, $i=1, \ldots,n$, we have $(\xi_{1}z_{1},\ldots,\xi_{n}z_{n}) \in \Omega$. Let $\mathcal{F}(\Omega)$ be the space of all holomorphic functions $f$ in $\Omega$ into $\mathbb{C}$. We write $\ell^n_{p}$ for the Banach space defined by $\mathbb{C}^n$ endowed with the $p$-norm $\norm{z}_{p}:= \left(\sum_{i=1}^{n}|z_{i}|^p\right)^{1/p},\, 1\leq q<\infty$ and $\norm{z}_{\infty}:=\sup^n_{i=1}|z_i|$. For $q \in [1,\infty]$, consider the unit balls in Minkowski space $\ell^n_{q}$ as
$$B_{\ell^n _q}= \left\{z \in \mathbb{C}^n:\norm{z}_{q}=\left(\sum_{i=1}^{n}|z_{i}|^q\right)^{1/q}<1\right\}\,\,\, \mbox{for}\,\,\,1\leq q<\infty$$ 
and $B_{\ell^n _\infty}= \{z \in \mathbb{C}^n: \norm{z}_{\infty}=\sup _{1\leq i \leq n} |z_{i}|<1\}$ which are Reinhardt domains of special interest in our context. For each $f\in \mathcal{F}(\Omega)$, we have the monomial series expansion of $f$ as $f(z)=\sum_{\alpha\in \mathbb{N}^n_{0}}x_{\alpha}(f)z^{\alpha}$ where the $\alpha$-th coefficient of the expansion given by  $x_{\alpha}(f)=(\partial^{\alpha}f(0)/\alpha!)$.
In 1997, Boas and Khavinson \cite{boas-1997} introduced the following multidimensional Bohr radius.

\begin{defn}
	The Bohr radius $K^n(\Omega)$ of a Reinhardt domain $\Omega\subset \mathbb{C}^n$  with respect to $\mathcal{F}(\Omega)$ is defined by
	\begin{equation} \label{h-e-1-a}
	K^n(\Omega):=\sup\left\{r\geq 0: \sup_{z\in r\Omega}\sum_{\alpha\in \mathbb{N}^n_{0}} \left|x_{\alpha}(f) z^{\alpha}\right| \leq  \norm{f}_{\Omega}\,\, \mbox{for all}\, f\in \mathcal{F}(\Omega)\right\},
	\end{equation}
where $\norm{f}_{\Omega}=\sup \{|f(z)|: z \in \Omega\}$ and $\mathbb{N}_0=\mathbb{N}\cup\{0\}$.
\end{defn}
\noindent
 We write $K^n(\Omega)=K(\Omega)$ for $n=1$. The celebrated theorem of Bohr \cite{Bohr-1914} states that $K(\mathbb{D})=1/3$. Although Bohr has obtained the radius $1/6$, which was later improved to $1/3$ by M. Riesz, I. Schur, and N. Wiener \cite{sidon-1927,tomic-1962}, independently. We usually say the inequality \eqref{h-e-1-a} is a Bohr inequality and the occurrence of this type inequality for all functions in $\mathcal{F}(\Omega)$ is known as Bohr phenomenon. When $\Omega=\mathbb{D}$, \eqref{h-e-1-a} is the classical Bohr inequality and $K(\mathbb{D})=1/3$ is the classical Bohr radius. Surprisingly, the exact value of the constant $K^n(\Omega)$ is not known for any other domain. The primary results of Boas and Khavinson \cite{boas-1997} and Boas \cite{boas-2000} have been able to provide a partial successful estimates for the Bohr radius $K^n(\Omega)$ for $\Omega=B_{\ell^n _q}$, $q \in [1,\infty]$. Their way of approaches towards finding the extimates of $K^n(B_{\ell^n _q})$ shows the difficulties in obtaining exact value of $K^n(B_{\ell^n _q})$. Therefore, it is always challenging to work on finding estimates of $K^n(\Omega)$ for an arbitrary Reinhardt domain. 
\vspace{1mm}

In this paper, we extensively study the arithmetic Bohr radius for certain class of holomorphic functions defined on complete Reinhardt domains in $\mathbb{C}^n$ which are not necessarily bounded. The notion of arithmetic Bohr radius has been first introduced by Defant {\it et. al.} \cite{defant-2007}. Let $\mathcal{B}(\Omega)$ be the class of all holomorphic functions $f: \Omega \rightarrow \mathbb{C}$ on a complete Reinhardt domain $\Omega\subset \mathbb{C}^n$ such that $\mbox{Re} (f(z))>0$ and $f(0)=1$. To study the arithmetic Bohr radius for the class $\mathcal{B}(\Omega)$, we introduce the following notion which will be the center of discussion in our paper.

\begin{defn}
	The	arithmetic Bohr radius of $\Omega$ with respect to the class $\mathcal{B}(\Omega)$, denoted by $A_p(\mathcal{B}(\Omega))$ and defined by
	\begin{equation*}
		A_{p}(\mathcal{B}(\Omega)):= \sup \left\{\frac{1}{n}\sum_{j=1}^{n}r_j\,|\, r\in \mathbb{R}^n_{\geq0}, \, \forall \, f\in \mathcal{B}(\Omega) : \frac{1}{2}\left(|c_0(f)|^p + \sum_{m=1}^{\infty}\sum_{|\alpha|=m}|c_{\alpha}r^{\alpha}|^p\right)^{\frac{1}{p}}\leq 1 \right\},
	\end{equation*}
	where $1\leq p< \infty$ and $\mathbb{R}^n_{\geq 0}=\{r=(r_1,.\,.\,.\,,r_n)\in \mathbb{R}^n : r_{i}\geq 0, 1\leq i\leq n\}$. We write $A_p(\Omega)$ for $A_p(\mathcal{B}(\Omega))$.
\end{defn}
It is worth to note that $A_p(\cdot)$ is increasing, that is, $A_p(\Omega_{1})\leq A_p(\Omega_{2})$ whenever $\Omega_{1} \subset \Omega_{2}.$ Let $\mathcal{P}(\Omega)$ be the set of all polynomials in $\mathcal{B}(\Omega)$ and $\mathcal{P}^m(\Omega)$ denote the set of all $m$-homogeneous polynomials in $\mathcal{B}(\Omega)$ defined on $\Omega$.
Arithmetic Bohr radius has rich properties and one of them is in describing the domain of existence of the monomial expansion of bounded holomorphic functions in a complete Reinhardt domain (see \cite{prengel-2005}). As an extended study recently, Kumar \cite{kumar-2023} has studied the arithmetic Bohr radius and answered certain questions raised by Defant {\it et. al.} in \cite{defant-2007}. 
\vspace{2mm}

The systematic and groundbreaking progress on Bohr problem for bounded holomorphic functions inspires us to study Bohr phenomenon problem for functions that are not necessarily bounded, more precisely for function whose images lie in the right half-plane. In 2000, Aizenberg {\it et al.} \cite{aizn-2000b} proved that if $f(z)=\sum_{k=0}^{\infty}a_{k}z^k$ be any holomorphic function with positive real part and $f(0)>0$, then
\begin{equation} \label{h-e-1}
	\sum_{k=0}^{\infty} |a_{k}z^k| \leq 2 f(0)	
\end{equation}
for $|z|\leq 1/3$ and the constant $1/3$ cannot be improved. It is worth mentioning that without loss of generality we can assume $f(0)=1$.  Later this work and \eqref{h-e-1} have been extended to several variable settings by Aizenberg {\it et al.} \cite{aizenberg-2003} while $p$-Bohr radius settings for functions in $\mathcal{B}(\Omega)$ in single variable settings have been extensively studied in \cite{Djakov & Ramanujan & J. Anal & 2000}. Here we consider \eqref{h-e-1} for functions in $\mathcal{B}(\Omega)$, where $\Omega$ is arbitrary Reinhardt domain in $\mathbb{C}^n$ and introduce the notion of Bohr radius for $\mathcal{B}(\Omega)$. 
\begin{defn}
	For $p>0$, the Bohr radius $H^{n}_{p}(\Omega)$ of $\Omega$ with respect to $\mathcal{B}(\Omega)$ defined by 
	\begin{equation}
		H^{n}_{p}(\Omega):=\sup\left\{r\geq 0:	\sup_{z\in r\Omega} \left(\frac{1}{2}\left(|c_0(f)|^p + \sum_{m=1}^{\infty}\sum_{|\alpha|=m}|c_{\alpha}(f)z^{\alpha}|^p\right)^{\frac{1}{p}}\right)\leq 1\right\},
	\end{equation}
	where $f \in \mathcal{B}(\Omega)$ with $f(z)=\sum_{\alpha\in \mathbb{N}^n_{0}} c_{\alpha}(f) z^{\alpha}$.
\end{defn}
Note that $H^{1}_{1}(\mathbb{D})=1/3$ while $H^{1}_{p}(\mathbb{D})=((2^p -1)/(2^{p+1}-1))^{1/p}$ for any $p>0$ (see \cite{Djakov & Ramanujan & J. Anal & 2000}). Motivated by the approaches in \cite{aizenberg-2003} and \cite{Djakov & Ramanujan & J. Anal & 2000}, Das \cite{das-2023} has recently considered \eqref{h-e-1} in more general setting for holomorphic functions in $B_{\ell^n_{\infty}}$ with positive real part. 
 For different values of $p$, $H^{n}_{p}(B_{\ell^n _\infty})$ has the following surprising asymptotic behavior due to \cite{das-2023}.
\begin{thm} \label{das-thm-2023} \cite{das-2023}
	For any $n>1$, 
	\begin{equation*}
		H^{n}_{p}(B_{\ell^n _\infty})= \left(\frac{2^p -1}{2^{p+1}-1}\right)^{\frac{1}{p}}
	\end{equation*}
	for $p \in [2,\infty)$ and 
	\begin{equation*}
		H^{n}_{p}(B_{\ell^n _\infty}) \sim \left(\frac{\log \, n}{n}\right)^{\frac{2-p}{2p}}
	\end{equation*}
	for $p \in (0,2)$.
\end{thm} 

The main aim of this paper is to find the exact value of $H^{n}_{p}(\Omega)$ in terms of arithmetic Bohr radius. To the best of our knowledge, nothing has been done to describe $H^{n}_{p}(\Omega)$ in terms of arithmetic Bohr radius.\\[2mm]

In the recent year, there has been a great progress in finding the exact value of multdimension Bohr radii. The problem of finding Bohr radius has been appeared in  several different contexts of mathematics, for instance, for Banach algebras and uniform algebras (see  \cite{paulsen-2002,paulsen-2004}), for complex manifolds (see \cite{aizn-2000b,aizenberg-2001a}), for ordinary and vector valued Dirichlet series (see \cite{bala-2006,defant-2008}), for elliptic equations (see \cite{aizenberg-2001b}), for Faber-Green condenser (see \cite{lassère-2017}), for free holomorphic functions (see \cite{popescu-2019}), for vector-valued holomorphic functions (see \cite{defant-2012}), for local Banach space theory (see \cite{defant-2003}), for domain of monomial convergence (see \cite{defant-2009}), for harmonic and pluriharmonic mappings (see \cite{hamada-2022-JFA}), Hardy spaces (see \cite{bene-2004}), and also in multidimensional settings (see \cite{aizn-2000a,boas-1997,boas-2000,defant-2004, defant-2006,Liu-Pon-PAMS-2020}).
The classical Bohr inequality was overlooked and did not get much attention for many years until it was used by Dixon \cite{Dixon & BLMS & 1995} to answer a long-standing open question related to Banach algebra satisfying a von Neumann inequality. In $1989$, Dineen and Timoney \cite{Dineen-Timoney-1989} first initiated the study of the constant $K^n(B_{\ell^n _\infty})$ and their result has been clarified by Boas and Khavinson in \cite{boas-1997}. In $1997$, Boas and Khavinson \cite{boas-1997} obtained the following lower and upper bounds of $K^n(B_{\ell^n _\infty})$ for each $n \in \mathbb{N}$ with $n \geq 2$,
\begin{equation}\label{Pal-Vasu-P3-e-1.1}
	\frac{1}{3\sqrt{n}} < K^n(B_{\ell^n _\infty}) < 2 \sqrt{\frac{\log n}{n}}.	
\end{equation}
The exact value of $K^n(B_{\ell^n _\infty})$ is still an open problem and the paper of Boas and Khavinson \cite{boas-2000} has aroused new interest in the multidimensional Bohr radius problem, and it has been a source of inspiration for many researchers to work further on this problem. Later, Aizenberg \cite{aizn-2000a} has obtained the following estimates of the constant $K^n(B_{\ell^n _1})$,
\begin{equation}\label{Pal-Vasu-P3-e-1.2}
	\frac{1}{3e^{1/3}} < K^n(B_{\ell^n _1}) \leq \frac{1}{3}.
\end{equation}
In $2000$, Boas \cite{boas-2000} extended the estimates \eqref{Pal-Vasu-P3-e-1.1} and \eqref{Pal-Vasu-P3-e-1.2} to $K^n(B_{\ell^n _q})$ for $1 <q<\infty$. For a fixed $n>1$, Boas \cite{boas-2000} has shown that, if $1\leq q <2$, then 
\begin{equation} \label{itdn-e-1}
	\frac{1}{3\sqrt[3]{e}}\left(\frac{1}{n}\right)^{1-\frac{1}{q}} \leq K^n(B_{\ell^n _q}) < 3 \left(\frac{\log \, n}{n}\right)^{1-\frac{1}{q}}
\end{equation}
and if $2 \leq q \leq \infty$, then
\begin{equation} \label{itdn-e-2}
	\frac{1}{3} \sqrt{\frac{1}{n}} \leq K^n(B_{\ell^n _q}) < 2 \sqrt{\frac{\log\,n}{n}}.
\end{equation}
In view of \eqref{itdn-e-1} and \eqref{itdn-e-2}, we see that the upper bounds contain a logarithmic factor but the lower bounds do not. For almost nine years, it was understood that the lower bound of \eqref{Pal-Vasu-P3-e-1.1}, \eqref{itdn-e-1}, and \eqref{itdn-e-2} could not be improved. Later, in $2006$, Defant and Frerick \cite{defant-2006} obtained a logarithmic lower bound which is almost correct asymptotic estimates for the Bohr radius $K^n(B_{\ell^n _q})$ with $1\leq q \leq \infty$. In particular, Defant and Frerick have proved that, if $1\leq q \leq \infty$ then there is a constant $c>0$ such that 
\begin{equation} \label{itdn-e-3}
	\frac{1}{c} \left(\frac{\log\,n/\log\, \log\,n}{n}\right)^{1-\frac{1}{\min(q,2)}} \leq K^n(B_{\ell^n _q}) \,\,\,\,\, \mbox{for all} \,\,\, n>1.
\end{equation}

\vspace{2mm}

From the above discussions, it seems that determining the exact value of multidimensional Bohr radius is very complicated even for the dimension $n=2$. As a part of our study, we give an independent way to investigate the multidimensional Bohr radius in terms of arithmetic Bohr radius for the class $\mathcal{B}(\Omega)$. In that process, we show that the arithmetic Bohr radius is closely related to classical Bohr radius.\\[1mm]

The organization of this paper is as follows. After the Introduction part, Section \ref{section-2} contains all the results of this paper and Section \ref{section-03} consists the proof of those results. We find exact value of $H^n_p(\Omega)$ for the unit ball in Minkowski space $\ell^n_p$ by establishing the relation between Bohr radius and arithmetic Bohr radius in Theorem \ref{Pal-Vasu-P4-thm-2.2}. We study the connection between arithmetic Bohr radius of $\mathcal{B}(\Omega)$ and arithmetic Bohr radius for the homogeneous polynomials in $\mathcal{B}(\Omega)$ in Proposition \ref{Pal-Vasu-P4-prop-2.1}. Moreover, we give estimates for arithmetic Bohr radius $A_p(B_{\ell^n_q})$ in terms of Bohr radius $H^n_p(B_{\ell^n_q})$ for $1\leq q<\infty$ and $q=\infty$ in Theorem \ref{Pal-Vasu-P4-thm-2.3} and Theorem \ref{Pal-Vasu-P4-thm-2.5} respectively. As a consequence, Corollary \ref{Pal-Vasu-P4-cor-2.5} shows that Bohr radius $H^1_1(B_{\ell^n_1})$ for unit ball of $\ell^n_1$ coincides with $1/3$, which is the classical Bohr radius.
\vspace{1mm}

\section{Main Results}\label{section-2}
In our first result, we provide an estimate for arithmetic Bohr radius of $\mathcal{B}(\Omega)$ in terms of the arithmetic Bohr radius for $m$-homogeneous polynomials in $\mathcal{B}(\Omega)$, where $\Omega$ being complete Reinhardt domain. 
\begin{prop}\label{Pal-Vasu-P4-prop-2.1}
	Let $\Omega$ be a complete Reinhardt domain in $\mathbb{C}^n$ and $1\leq p <\infty$. Then we have 
	\begin{equation}\label{Pal-Vasu-e-2.2}
		 \frac{1}{3^{1/p}}\, A_p\left(\bigcup_{m=1}^{\infty}\mathcal{P}^{m}(\Omega)\right) \leq A_p\left(\mathcal{B}(\Omega)\right) \leq A_p\left(\bigcup_{m=1}^{\infty}\mathcal{P}^m(\Omega)\right).
	\end{equation}  
\end{prop}

We present the next main result as Theorem \ref{Pal-Vasu-P4-thm-2.2} where we obtain the exact value of $n$-dimensional Bohr radius $H^n_{p}(B_{\ell^n_{q}})$ in terms of the arithmetic Bohr radius $A_p(B_{\ell^n_{q}})$ for the unit ball in $\ell^n_{q}$-spaces. Before briefing Theorem \ref{Pal-Vasu-P4-thm-2.2}, we establish a relation between the arithmetic Bohr radius $A_p(\Omega)$ and the Bohr radius $H^n_{p}(\Omega)$ for bounded Reinhardt domain $\Omega$ in $\mathbb{C}^n$, which we offer as Lemma \ref{Pal-Vasu-P4-lem-2.1}. To make the statement precise, we require the following notation from \cite{defant-2004}. For  bounded Reinhardt domains $\Omega_{1}, \Omega_{2}\subset \mathbb{C}^n$, let
\begin{equation*}
	S(\Omega_1, \Omega_{2}):= \inf \left\{t>0 : \Omega_{1} \subset t\Omega_{2}\right\}.
\end{equation*}
 By a Banach sequence space $X$, we mean a complex Banach space $X\subset \mathbb{C}^{\mathbb{N}}$ such that $\ell_{1}\subset X \subset \ell_{\infty}$. If $\Omega$ is a bounded Reinhardt domain in $\mathbb{C}^n$ and $X$ and $Y$ are Banach sequence spaces we write
\begin{equation}\label{Pal-Vasu-P4-e-2.3}
	S(\Omega, B_{X_n})= \sup_{z\in \Omega}\norm{z}_{X} \quad \mbox{and }\quad  S(B_{X_n}, B_{Y_n})= \norm{\mbox{id}: X_n \rightarrow Y_n},
\end{equation}
where $X_n$(resp. $Y_n$) is the space spanned by first $n$ canonical basis vectors $e_n$ in $X$(resp. $Y$).

\begin{rem}
	For a bounded Reinhardt domain $\Omega$ in $\mathbb{C}^n$, observe that $S(\Omega, t\Omega)=1/t$ and $S(t\Omega, \Omega)=t$  for all $t>0.$ 

\end{rem}
The following lemma relates the Bohr radius $H^n_{p}(\Omega)$ and the arithmetic Bohr radius $A_p(\Omega)$ for bounded Reinhardt domain $\Omega$.
\begin{lem}\label{Pal-Vasu-P4-lem-2.1}
	Let $\Omega \subset\mathbb{C}^n$ be a bounded Reinhardt domain in $\mathbb{C}^n$ and $1\leq p<\infty$. Then we have 
	\begin{equation*}
		A_p(\Omega) \geq \frac{S(\Omega, B_{\ell^n_1})}{n}H^n_{p}(\Omega).
	\end{equation*}
\end{lem}
As discussed before, now we show the exact value of Bohr radius $H^n_{p}(\Omega)$ in terms of the arithmetic Bohr radius $A_p(\Omega)$ for $\Omega=B_{\ell^n_{q}}$, $1\leq q\leq \infty$.
\begin{thm}\label{Pal-Vasu-P4-thm-2.2}
	Let $1\leq  p< \infty$. Then for every $1\leq q\leq \infty$ and for all $n\in \mathbb{N}$, we have 
	\begin{equation*}
		A_p(B_{\ell^n_q})= \frac{H^n_p(B_{\ell^n_q})}{n^{1/q}}.
	\end{equation*}
\end{thm}
Next, we obtain an interesting relation between the classical Bohr radius $H^1_{p}(\mathbb{D})$ and the arithmetic Bohr radius $A_p(B_{\ell^n_{q}})$ for $1\leq q<\infty$. Further, we shall see that this relation helps us to compare the classical Bohr radii for unit disk and unit ball in $\mathbb{C}^n$. 
\begin{thm}\label{Pal-Vasu-P4-thm-2.3}
Let $1\leq p<\infty$. Then for every $n\in \mathbb{N}$ and $1\leq q<\infty$ we have 
\begin{equation*}
	\frac{H^1_p(\mathbb{D})}{n}\leq A_p\left(B_{\ell^n_q}\right)\leq \left(\frac{H^1_p(\mathbb{D})}{n^{1/p}}\right)^{1/q}.
\end{equation*}
\end{thm}
In view of Theorem \ref{Pal-Vasu-P4-thm-2.2} and Theorem \ref{Pal-Vasu-P4-thm-2.3}, we obtain the following interesting estimate.
\begin{thm}\label{Pal-Vasu-P4-thm-2.4}
	For every $1\leq p,q <\infty$ and $n \in \mathbb{N}$, we have 
	\begin{equation*}
		\frac{H^1_{p}(\mathbb{D})}{n^{1-(1/q)}}\leq H^n_{p}(B_{\ell^n_{q}})\leq \left(\frac{H^1_{p}(\mathbb{D})}{n^{(1/p)}-1}\right)^{1/q}.
	\end{equation*}
\end{thm}
The exact value of Bohr radius $H^n_{p}(B_{\ell^n_{\infty}})$ for the unit polydisc has been studied by Das \cite{das-2023} as we have seen in Theorem \ref{das-thm-2023}, whereas the exact value for unit polyballs in $\ell^n_{q}$-spaces $(1\leq q<\infty)$ is still an open problem. In view of Theorem \ref{Pal-Vasu-P4-thm-2.4}, we observe that the exact value of Bohr radius $H^1_{1}(B_{\ell^n_{1}})$ for the unit ball in $\ell^n_{1}$ space is exactly $1/3$.
\begin{cor}\label{Pal-Vasu-P4-cor-2.5}
	For every $n\in \mathbb{N}$, we have 
	\begin{equation*}
		H^1_{1}(B_{\ell^n_{1}})=H^1_1(\mathbb{D})=1/3.
	\end{equation*}
\end{cor}
We also study the case $q=\infty$ in Theorem \ref{Pal-Vasu-P4-thm-2.3} and obtain the following estimate for the arithmetic Bohr radius $A_p(B_{\ell^n_{\infty}})$ in terms of the classical Bohr radius $H^1_{p}(\mathbb{D})$.
\begin{thm}\label{Pal-Vasu-P4-thm-2.5}
Let $1\leq p <\infty$. Then for each $n\in \mathbb{N}$, we have   
	\begin{equation*}
		\frac{H^1_p(\mathbb{D})}{n} \leq A_p(B_{\ell^n_{\infty}}) \leq \frac{H^1_p(\mathbb{D})}{n^{(1/p)-1}}.
	\end{equation*}
\end{thm}
In the following section, we present the proof of Proposition \ref{Pal-Vasu-P4-prop-2.1}, Lemma \ref{Pal-Vasu-P4-lem-2.1}, Theorem \ref{Pal-Vasu-P4-thm-2.2}, Theorem \ref{Pal-Vasu-P4-thm-2.3} and Theorem \ref{Pal-Vasu-P4-thm-2.5}.

\section{Proof of Main Results}\label{section-03}
\begin{pf} [{\bf Proof of Proposition \ref{Pal-Vasu-P4-prop-2.1}}]
		Since we have the following inclusion 
	\begin{equation*}
		\bigcup_{m=1}^{\infty} \mathcal{P}^{m}(\Omega) \subset\mathcal{B}(\Omega),
	\end{equation*}
the right-hand inequality of \eqref{Pal-Vasu-e-2.2},
\begin{equation}
	A_p\left(\mathcal{B}(\Omega)\right) \leq A_p\left(\bigcup_{m=1}^{\infty}\mathcal{P}^m(\Omega)\right)
\end{equation}
 holds.
	Choose $r\in \mathbb{R}^n_{\geq 0}$ be such that for all $m$-homogeneous polynomial $g_m \in \mathcal{P}^m(\mathbb{C}^n)$ contained in $\mathcal{B}(\mathbb{C}^n)$,
	\begin{equation}\label{Pal-Vasu-P4-e-2.1}
		\frac{1}{2}\left(\sum_{|\alpha|=m}|c_{\alpha}(g_m)|^pr^{p\alpha}\right)^{\frac{1}{p}}\leq 1.
	\end{equation}
	Our aim is to show that 
	\begin{equation}\label{Pal-Vasu-P4-e-2.2}
		\frac{1}{3^{1/p}}\,\sum_{i=1}^{n}r_i \leq A_p(\mathcal{B}(\Omega)).
	\end{equation}
	Take $f(z)=\sum_{\alpha\in \mathbb{N}^n_{0}}c_{\alpha}(f)z^{\alpha} \in \mathcal{B}(\Omega)$. Then, in view of \eqref{Pal-Vasu-P4-e-2.1}, we obtain
	\begin{align*}
		\frac{1}{2}\left(\sum_{\alpha\in \mathbb{N}^n_{0}}|c_{\alpha}(f)|^p\left(\frac{r^p}{3}\right)^{\alpha}\right)^{1/p} &= \frac{1}{2} \left(|c_0(f)|^p + \sum_{m=1}^{\infty}\sum_{|\alpha|=m}|c_{\alpha}(f)|^p\left(\frac{r^p}{3}\right)^{\alpha}\right)^{1/p}\\ & = \frac{1}{2}\left(|c_0(f)|^p + \sum_{m=1}^{\infty}\frac{1}{3^m}\sum_{|\alpha|=m}|c_{\alpha}(f)r^{\alpha}|^p\right)^{1/p}\\ & \leq \frac{1}{2}\left(1+ 2^p\sum_{m=1}^{\infty}\frac{1}{3^m}\right)^{1/p}=\frac{1}{2}\left(1+2^{p-1}\right)^{1/p}\leq 1,
	\end{align*}
	which gives the estimate \eqref{Pal-Vasu-P4-e-2.2}. Hence, 
	\begin{equation*}
		\frac{1}{3^{1/p}}A_p(\mathcal{P}^m(\Omega))\leq A_p(\mathcal{B}(\Omega)) \quad \mbox{for all}\,\,\, m\geq 1.
	\end{equation*}
	As a consequence, we obtain the left-hand inequality of \eqref{Pal-Vasu-e-2.2}. This completes the proof.
\end{pf}

\begin{pf} [{\bf Proof of Lemma \ref{Pal-Vasu-P4-lem-2.1}}]
	By the virtue of \eqref{Pal-Vasu-P4-e-2.3}, we have 
	\begin{equation*}
		S(\Omega, B_{\ell^n_{1}})=\sup_{z\in \Omega}\norm{z}_{\ell^n_{1}}.
	\end{equation*}
	Thus for given $0<\epsilon<H^n_{p}(\Omega)$, we can find an element $z_0\in \Omega$ such that 
	\begin{equation*}
		\norm{z_0}_{\ell^n_1}\geq S(\Omega, B_{\ell^n_{1}})-\epsilon.
	\end{equation*}
	Let $t:=H^n_{p}(\Omega)-\epsilon$, $v:=sz_0$, and $r:=s|z_0|=|v|.$ Since $v\in t\Omega$ and $t<H^n_p(\Omega)$, for $f=\sum_{\alpha\in \mathbb{N}^n_{0}}c_{\alpha}(f)z^{\alpha}\in \mathcal{B}(\Omega)$, we have
	\begin{equation*}
		\frac{1}{2}\left(|c_0(f)|^p + \sum_{m=1}^{\infty}\sum_{|\alpha|=m}|c_{\alpha}(f)|^pr^{p\alpha}\right)^{1/p}=\frac{1}{2}\left(|c_0(f)|^p + \sum_{m=1}^{\infty}\sum_{|\alpha|=m}|c_{\alpha}v^{\alpha}|^p\right)^{1/p}\leq 1.
	\end{equation*}
	Therefore, we obtain 
	\begin{equation*}
		A_p(\Omega) \geq \frac{1}{n}\sum_{i=1}^{n}r_i = \frac{\norm{r}_{1}}{n} \frac{H^n_p(\Omega)-\epsilon}{n}\norm{z_0}_{\ell^n_{1}}\geq \frac{H^n_p(\Omega)-\epsilon}{n}\left(S(\Omega, B_{\ell^n_{1}})-\epsilon\right)
	\end{equation*}
	holds for all $\epsilon>0.$ Letting $\epsilon \rightarrow 0$, we have 
	\begin{equation*}
		A_p(\Omega) \geq \frac{S(\Omega, B_{\ell^n_1})}{n}H^n_{p}(\Omega).
	\end{equation*}
	This completes the proof.
\end{pf}

\begin{pf} [{\bf Proof of Theorem \ref{Pal-Vasu-P4-thm-2.2}}]
In view of H\"{o}lder's inequality, we have	$S(B_{\ell^n_q}, B_{\ell^n_1})=n^{1-(1/q)}$. Using this fact in Lemma \ref{Pal-Vasu-P4-lem-2.1}, we obtain the inequality
	\begin{equation*}
		A_p(B_{\ell^n_q}) \geq  \frac{H^n_p(B_{\ell^n_q})}{n^{1/q}}.
	\end{equation*}
	Therefore, we need to show that
	\begin{equation}\label{Pal-Vasu-P4-e-2.5}
		A_p(B_{\ell^n_q}) \leq  \frac{H^n_p(B_{\ell^n_q})}{n^{1/q}}.
	\end{equation}
	Let $r=(r_1,.\,.\,., r_n)\in \mathbb{R}^n_{\geq 0}$ be such that for all $h(z)=\sum_{\alpha\in \mathbb{N}^n_{0}}c_{\alpha
	}(h)z^{\alpha}\in\mathcal{B}(B_{\ell^n_q})$, 
	\begin{equation*}
		\frac{1}{2}\left(|c_0(h)|^p + \sum_{m=1}^{\infty}\sum_{|\alpha|=m}|c_{\alpha}(h)|^pr^{p\alpha}\right)^{1/p}\leq 1.
	\end{equation*}
	To prove \eqref{Pal-Vasu-P4-e-2.5}, it suffices to prove that
	\begin{equation*}
		n^{\frac{1}{q}-1}\norm{r}_{1}\leq H^n_p(B_{\ell^n_q}).
	\end{equation*}
	Let $f\in \mathcal{B}(B_{\ell^n_q})$. It is worth noting that for $u\in n^{(1/q)-1}\norm{r}_1\overline{B_{\ell^n_q}}$, we have $\norm{u}_1\leq \norm{r}_1$.
	Therefore,
	\begin{equation*}
		\frac{1}{2}\left(|c_0(f)|^p + \sum_{m=1}^{\infty}\sum_{|\alpha|=m}|c_{\alpha}(f)u^{\alpha}|^p\right)^{1/p}\leq 1
	\end{equation*}
	for every $u \in n^{(1/q)-1}\norm{r}_1\overline{B_{\ell^n_q}}.$ So, we obtain $n^{(1/q)-1}\norm{r}_{1}\leq H^n_{p}(B_{\ell^n_{q}})$. Consequently, it follows that
	\begin{equation*}
		n^{\frac{1}{q}}A_p(B_{\ell^n_q})\leq H^n_p(B_{\ell^n_q}),
	\end{equation*}
	which gives our conclusion. This completes the proof.
\end{pf}

\begin{pf} [{\bf Proof of Theorem \ref{Pal-Vasu-P4-thm-2.3}}]
	First we show the left-hand inequality
	\begin{equation*}
		\frac{H^1_p(\mathbb{D})}{n}\leq A_p\left(B_{\ell^n_q}\right).
	\end{equation*}
	Assume $r=H^1_p(\mathbb{D})$ and $f\in \mathcal{B}(B_{\ell^n_q})$. We define $g(z)=f(ze_1)=f(z,0,.\,.\,.,0)$ for $z\in \mathbb{D}.$ Then $g:\mathbb{D}\rightarrow\mathbb{C}$ will be a holomorphic function on $\mathbb{D}$ with $\mbox{Re}(g(z))>0$ and $g(0)=1$. Therefore,
	\begin{align*}
		\frac{1}{2}\left\{|c_0(f)|^p+\sum_{k=1}^{\infty}\sum_{|\alpha|=k}|c_{\alpha}(f)|(r,0,.\,.\,.\,,0)^{p\alpha}\right\}^{\frac{1}{p}}= \frac{1}{2}\left\{|c_0(g)|^p+\sum_{k=1}^{\infty}|c_{k}(g)|r^{pk}\right\}^{\frac{1}{p}}\leq 1
	\end{align*}
	for all $f(z)=\sum_{\alpha\in \mathbb{N}^n_{0}}c_{\alpha}(f)z^{\alpha} \in \mathcal{B}(B_{\ell^n_q})$. Hence, we obtain $r/n\leq A_p(B_{\ell^n_q})$, which gives our desired inequality.\\[2mm]
	On the other hand, we want to prove that 
	\begin{equation*}
		A_p\left(B_{\ell^n_q}\right)\leq \left(\frac{H^1_p(\mathbb{D})}{n^{1/p}}\right)^{1/q}.
	\end{equation*}
	Let $r\in \mathbb{R}^n_{\geq 0}$ be such that for all $u\in \mathcal{B}(B_{\ell^n_q})$, we have
	\begin{equation*}
		\frac{1}{2}\left\{|c_0(u)|^p + \sum_{k=1}^{\infty}\sum_{|\alpha|=k}|c_{\alpha}(u)|^pr^{p\alpha}\right\}^{\frac{1}{p}}\leq 1.
	\end{equation*}
	Now it is sufficient to show that
	\begin{equation}\label{Pal-Vasu-P4-e-3.5}
		\frac{1}{n}\left(\sum_{j=1}^{n}r_j\right) \leq  \left(\frac{H^1_p(\mathbb{D})}{n^{1/p}}\right)^{1/q}.
	\end{equation}
	Fix $f\in \mathcal{B}(\mathbb{D})$, and define the function
	\begin{equation*}
		v(z)=z^q_1+\cdots+z^q_n, \quad z=(z_1,.\,.\,.,z_n)\in B_{\ell^n_q}. 
	\end{equation*}
	 Consider $u=f\circ v$ such that $\mbox{Re}(u(z))=\mbox{Re}(f(v(z)))>0$ and $u(0)=f(v(0))=f(0)=1$. Moreover, for each $z\in B_{\ell^n_q}$, we have 
	\begin{equation*}
		u(z)=\sum_{k=0}^{\infty}c_{k}(f)v(z)^k=\sum_{k=0}^{\infty}c_k(f)\sum_{|\alpha|=k}\frac{k!}{\alpha!}z^{q\alpha}=\sum_{\alpha\in \mathbb{N}^n_{0}}c_{\alpha}(u)z^{q\alpha},
	\end{equation*}
	where $c_{\alpha}(u)=c_k(f)(k!/\alpha!)$ whenever $|\alpha|=k$. Then for all $z\in B_{\ell^n_q}$, we have
	\begin{align*}
		\frac{1}{2}\left\{|c_0(u)|^p + \sum_{k=1}^{\infty}\sum_{|\alpha|=k}|c_{\alpha}(u)z^{q\alpha}|^p\right\}^{\frac{1}{p}} &= \frac{1}{2}\left\{|c_0(f)|^p+\sum_{k=1}^{\infty}|c_{k}(f)|^p\sum_{|\alpha|=k}\left(\frac{k!}{\alpha!}\right)^p|z|^{pq\alpha}\right\}^{\frac{1}{p}}\\& \geq 
		\frac{1}{2}\left\{|c_0(f)|^p+\sum_{k=1}^{\infty}|c_{k}(f)|^p\sum_{|\alpha|=k}\frac{k!}{\alpha!}|z|^{pq\alpha}\right\}^{\frac{1}{p}}\\ & = 
		\frac{1}{2}\left\{|c_0(f)|^p+\sum_{k=1}^{\infty}|c_{k}(f)|^p\norm{z}^{pqk}_{pq}\right\}^{\frac{1}{p}}
	\end{align*}
	so that finally we have 
	\begin{align*}
		\frac{1}{2}\left\{|c_0(f)|^p+\sum_{k=1}^{\infty}|c_{k}(f)|^p\norm{r}^{pqk}_{pq}\right\}^{\frac{1}{p}}\leq \frac{1}{2}\left\{|c_0(u)|^p + \sum_{k=1}^{\infty}\sum_{|\alpha|=k}|c_{\alpha}(u)r^{q\alpha}|^p\right\}^{\frac{1}{p}}\leq 1. 
	\end{align*}
	It follows that $\norm{r}^q_{pq}\leq H^1_{p}(\mathbb{D})$. By virtue of H\"{o}lder's inequality, we have $\norm{r}^{pq}_{1}\leq n^{pq-1}\norm{r}^{pq}_{pq}$. Hence, we obtain $n^{1-pq}\norm{r}^{pq}_{1}\leq (H^1_{p}(\mathbb{D}))^p$, which gives the estimate \eqref{Pal-Vasu-P4-e-3.5}. This completes the proof. 
\end{pf}

\begin{pf} [{\bf Proof of Theorem \ref{Pal-Vasu-P4-thm-2.5}}]
	Let $r=H^1_{p}(\mathbb{D})$ and $f(z)=\sum_{\alpha\in \mathbb{N}^n_{0}}c_{\alpha}(f)z^{\alpha} \in \mathcal{B}(B_{\ell^n_{\infty}})$. Consider the function $g(z)=f(\xi z)$, where $\xi=(1, 0,.\,.\,.,0)$ and $z\in \mathbb{D}$. Clearly, $g$ is an holomorphic function on unit disk $\mathbb{D}$ with $\mbox{Re}(g(z))>0$ and $g(0)=1$. 
	Then we have 
	\begin{equation*}
		\frac{1}{2}\left(|c_0(f)|^p + \sum_{k=1}^{\infty}\sum_{|\alpha|=k}|c_{\alpha}(f)(r,0,.\,.\,.,0)^{\alpha}|^{p}\right)= \frac{1}{2}\left(|c_0(g)|^p + \sum_{k=1}^{\infty}|c_k(g)|^pr^{pk}\right)\leq 1.
	\end{equation*}
	Therefore, it gives us $(r/n)\leq A_p(B_{\ell^n_{\infty}})$, and hence we obtain $	(H^1_p(\mathbb{D})/n) \leq A_p(B_{\ell^n_{\infty}}).$ Conversely, we prove that 
	\begin{equation*}
		A_p(B_{\ell^n_{\infty}}) \leq \frac{H^1_p(\mathbb{D})}{n^{1/p-1}}.
	\end{equation*}
	Suppose $r\in \mathbb{R}^n_{\geq 0}$ such that for all $h \in \mathcal{B}$,
	\begin{equation*}
		\frac{1}{2}\left(|c_0(h)|^p + \sum_{k=1}^{\infty}\sum_{|\alpha|=k}|c_{\alpha}(h)|^pr^{p\alpha}\right)^{\frac{1}{p}} \leq 1.
	\end{equation*}
	Let $f:\mathbb{D}\rightarrow\mathbb{C}$ be a holomorphic function such that $\mbox{Re} f(z)>0$ and $f(0)=1$. Now we consider the function $s: B_{\ell^n_{\infty}}\rightarrow \mathbb{D}$ defined by 
	\begin{equation*}
		s(z)=\frac{1}{n}(z_1+\cdots+z_n),\quad z\in B_{\ell^n_{\infty}}.
	\end{equation*}
	Now if we set $h=f\circ s$, then we have $h\in \mathcal{B}(B_{\ell^n_{\infty}})$ with $\mbox{Re}(h(z))>0$ and $h(0)=1$. Also, for each $z\in B_{\ell^n_{\infty}}$,
	\begin{equation*}
		h(z)=\sum_{k=1}^{\infty}c_k(f)s(z)^k=\sum_{k=0}^{\infty}\frac{c_k(f)}{n^k}\sum_{|\alpha|=k}\frac{k!}{\alpha!}z^{\alpha}=\sum_{\alpha\in \mathbb{N}^n_{0}}c_{\alpha}(h)z^{\alpha},
	\end{equation*}
	where 
	\begin{equation*}
		c_{\alpha}(h)=\frac{k!}{\alpha!}\left(\frac{c_k(f)}{n^k}\right)
	\end{equation*}
	whenever $|\alpha|=k$. Then for all $z\in B_{\ell^n_{\infty}}$, we have
	\begin{align*}
		\frac{1}{2}\left(|c_0(h)|^p + \sum_{k=1}^{\infty}\sum_{|\alpha|=k}|c_{\alpha}(h)z^{\alpha}|^p\right)^{\frac{1}{p}} & = \frac{1}{2}\left(|c_0(f)|^p+ \sum_{k=1}^{\infty}\frac{|c_k(f)|^p}{n^{kp}}\sum_{|\alpha|=k}\left(\frac{k!}{\alpha!}\right)^{p} z^{p\alpha}\right)^{\frac{1}{p}}\\[2mm] & \geq 
		\frac{1}{2}\left(|c_0(f)|^p+ \sum_{k=1}^{\infty}\frac{|c_k(f)|^p}{n^{kp}}\sum_{|\alpha|=k}\left(\frac{k!}{\alpha!}\right)z^{p\alpha}\right)^{\frac{1}{p}}\\[2mm]& = 
		\frac{1}{2}\left(|c_0(f)|^p + \sum_{k=1}^{\infty}\frac{|c_k(f)|^p}{n^{pk}}\norm{z}^{pk}_{p}\right)^{\frac{1}{p}}.
	\end{align*}
	Finally, we observe that
	\begin{equation*}
		\frac{1}{2}\left(|c_0(f)|^p + \sum_{k=1}^{\infty}\frac{|c_k(f)|^p}{n^{pk}}\norm{r}^{pk}_{p}\right)^{\frac{1}{p}}\leq \frac{1}{2}\left(|c_0(h)|^p + \sum_{k=1}^{\infty}\sum_{|\alpha|=k}|c_{\alpha}(h)|^pr^{pk}\right)^{\frac{1}{p}}\leq 1.
	\end{equation*}
	This shows that $(1/n)\norm{r}_p\leq H^1_{p}(\mathbb{D}).$ Again we have $\norm{r}^p_{1}\leq n^{p-1}\norm{r}^p_{p}.$ Hence we obtain
	\begin{equation*}
		\frac{1}{n}\norm{r}_1\leq n^{1-(1/p)}H^1_{p}(\mathbb{D}),
	\end{equation*}
which gives our desired inequality. This completes the proof.
\end{pf}

\noindent\textbf{Acknowledgment:} 
The research of the first named author is supported by SERB-CRG (DST), Govt. of India.  The research of the second named author is supported by Institute Post-Doctoral Fellowship of IIT Bombay, and the research of the third named author is supported by DST-INSPIRE Fellowship (IF 190721),  New Delhi, India.


\begin{thebibliography}{99}
	
	
	
	
	
	
	
	
	
	
	
	
	\bibitem{aizn-2000a} {\sc L. Aizenberg}, Multidimensional analogues of Bohr's theorem on power series, \textit{Proc. Amer. Math. Soc.} {\bf 128} (2000), 1147--1155.
	
	\bibitem{aizn-2000b} {\sc L. Aizenberg, A. Aytuna}, and {\sc P. Djakov}, An abstract approach to Bohr's phenomenon, {\it Proc. Amer. Math. Soc.} {\bf 128} (2000), 2611--2619.
	
	
	\bibitem{aizenberg-2001a} {\sc L. Aizenberg, A. Aytuna},  and {\sc P. Djakov}, Generalization of theorem on Bohr for bases in spaces of holomorphic functions of several complex variables, 
	{\it J. Math. Anal. Appl.} {\bf  258} (2001), 429--447.
	
	\bibitem{aizenberg-2001b} {\sc L. Aizenberg}  and {\sc N. Tarkhanov}, A Bohr Phenomenon for elliptic equations, 
	{\it Proc. Lond. Math. Soc.} {\bf  82} (2001), 385--401.
	
	\bibitem{aizenberg-2003} {\sc L. Aizenberg, I. B. Grossman}  and {\sc Yu. F. Korobeinik}, Some remarks on the Bohr radius	for power series (Russian), {\it Izv. Vyssh. Uchebn. Zaved. Mat.}, 2002, no. 10, 3–10; {\it translation	in Russian Math. (Iz. VUZ)}, 46 (2002), no. 10, 1–8 (2003).
	
	
	
	
	
	
	
	
	
	
	
	
	
	
	
	
	

	
	
	
	
	
	
	
	
	
	
	
	
	
	\bibitem{bala-2006} {\sc R. Balasubramanian}, {\sc B. Calado} and {\sc H. Queffélec}, The Bohr inequality for ordinary Dirichlet series, {\it Studia Math.} {\bf 175} (2006), 285--304.
	
	\bibitem{bene-2004} {\sc C. B{\'e}n{\'e}teau}, {\sc A. Dahlner} and {\sc D. Khavinson}, Remarks on the Bohr phenomenon, {\it Comput. Methods Funct. Theory} \textbf{4}(1) (2004), 1-19.
	
	
	
	
	
	
	
	\bibitem{boas-1997} {\sc H. P. Boas} and {\sc D. Khavinson}, Bohr's power series theorem in several variables, {\it Proc. Amer. Math. Soc.}  {\bf 125} (1997), 2975--2979.
	
		\bibitem{boas-2000} {\sc H. P. Boas}, Majorant Series,  {\it J. Korean Math. Soc.}  {\bf 37} (2000), 321--337.
	
	\bibitem{Bohr-1914} {\sc H. Bohr}, A theorem concerning power series,  {\it Proc. Lond. Math. Soc.} s2-13 (1914), 1--5.
	
	
	
	
	
	
	
	
	\bibitem{hamada-2022-JFA} {\sc S. Chen} and {\sc  H. Hamada}, Some sharp Schwarz-Pick type estimates and their applications of harmonic and pluriharmonic functions, {\it J. Funct. Anal.} {\bf 282} (2022), 109254, 42 pp.
	
	\bibitem{das-2023} {\sc N. Das}, Estimates for generalized Bohr radii in one and higher dimensions, (2023) https://arxiv.org/pdf/2204.12706.pdf.
	
	\bibitem{defant-2003} {\sc A. Defant, D. Garc\'{i}a}, and {\sc  M. Maestre}, Bohr power series theorem and local Banach space theory, {\it J. Reine Angew. Math.} {\bf 557} (2003), 173--197.
	
	\bibitem{defant-2004} {\sc A. Defant, D. Garc\'{i}a}, and {\sc  M. Maestre}, Estimates for the first and second Bohr radii of Reinhardt domains, {\it J. Appr. Theory} {\bf 128} (2004), 53--68.
	
	\bibitem{defant-2006} {\sc A. Defant} and {\sc L. Frerick}, A logarithmic lower bound for multi-dimenional Bohr radii, {\it Israel J. Math.} {\bf 152} (2006), 17--28.
	
	\bibitem{defant-2007} {\sc A. Defant,  M. Maestre}, and {\sc  C. Prengel}, The arithmetic Bohr radius, {\it Q. J. Math.} {\bf 59} (2008), 189--205.
	
	\bibitem{defant-2008} {\sc A. Defant, D. Garc\'{i}a, M. Maestre}, and {\sc D. P\'{e}rez-Garc\'{i}a}, Bohr's strip for vector-valued Dirichlet series, {\it Math. Ann.} {\bf 342} (2008), 533--555.
	
	 \bibitem{defant-2009} {\sc A. Defant,  M. Maestre}, and {\sc  C. Prengel}, Domains of convergence for monomial expansions of holomorphic functions in infinitely many variables, {\it J. Reine Angew. Math.} {\bf 634} (2009), 13--49.
	
	
	\bibitem{defant-2012} {\sc A. Defant, M. Maestre}, and {\sc  U. Schwarting}, Bohr radii of vector valued holomorphic functions, {\it Adv. Math.} {\bf 231} (2012), 2837--2857.
	
	
	
	
	
	\bibitem{Dixon & BLMS & 1995} {\sc P. G. Dixon}, Banach algebras satisfying the non-unital von Neumann inequality, \textit{Bull. Lond. Math. Soc.} \textbf{27} (4) (1995), 359--362.
	
	\bibitem{Djakov & Ramanujan & J. Anal & 2000} {\sc P. B. Djakov} and {\sc M. S. Ramanujan}, A remark on Bohr's theorem and its generalizations, \textit{J. Anal.} \textbf{8} (2000), 65--77.
	
	\bibitem{Dineen-Timoney-1989} {\sc S. Dineen} and {\sc R. M. Timoney}, Absolute bases, tensor products and a theorem of Bohr, \textit{Studia Math.} \textbf{94} (1989), 227--234.
	
	
	
	
	
	
	
	
	
	
	
	
	
	
	
	
	
	
	
	
	
	
	
	
	
	
	
	
	
	
	
	
	
	
	
	
	
	
	\bibitem{kumar-2023} {\sc S. Kumar}, On the multidimensional Bohr radius, {\it Proc. Amer. Math. Soc.} {\bf 151} (2023), 2001--2009.
	
	\bibitem{lassère-2017} {\sc P. Lassère} and {E. Mazzilli}, Estimates for the Bohr Radius of a Faber–Green Condenser in the Complex Plane, {\it Constr. Approx.} {\bf 45} (2017), 409--426.
	
	
	
	
	
	
	
	
	
	
	\bibitem{Liu-Pon-PAMS-2020} {\sc M. S. Liu} and {\sc S. Ponnusamy}, Multidimensional analogues of refined Bohr's inequality, {\it Proc. Amer. Math. Soc.} {\bf 149} (2021), 2133--2146.
	
	
	
	
	
	
	
	\bibitem{paulsen-2002} {\sc V. I. Paulsen, G. Popescu}, and {\sc D. Singh}, On Bohr's inequality, {\it Proc. Lond. Math. Soc.} s3-85 (2002), 493--512.
	
	\bibitem{paulsen-2004} {\sc V. I. Paulsen} and {\sc D. Singh}, Bohr’s inequality for uniform algebras, {\it Proc. Amer. Math. Soc.}, {\bf 132} (2004), 3577-3579.
	
	\bibitem{prengel-2005} {\sc C. Prengel}, {\it Domains of convergence in infinite dimensional holomorphy}, Ph.D. Thesis, University of Oldenburg, 2005.
	
	
	
	
	\bibitem{popescu-2019} {\sc G. Popescu}, Bohr inequalities for free holomorphic functions on polyballs, {\it Adv. Math.} {\bf 347} (2019), 1002-1053.
	
	
	
	
	
	
	
	
	
		\bibitem{sidon-1927} {\sc S. Sidon}, Uber einen satz von Hernn Bohr, {\it Math. Zeit.}  {\bf 26} (1927),  731--732.
	
	
	
	
	
	
	
	
		\bibitem{tomic-1962} {\sc M. Tomic}, Sur un theoreme de H. Bohr,  {\it Math. Scand.}  {\bf 11} (1962), 103--106.
	
	
	
	
\end{thebibliography}
\end{document}